\documentclass[12pt,leqno]{amsart} 

\usepackage{amssymb,amsmath} 
\usepackage{amscd}

\theoremstyle{plain} 
 \newtheorem{theorem}{\indent\sc Theorem}[section]
 
 \newtheorem{corollary}[theorem]{\indent\sc Corollary}
 \newtheorem{proposition}[theorem]{\indent\sc Proposition}

\theoremstyle{definition}

 \newtheorem{example}[theorem]{\indent\sc Example}

\newtheorem*{remark0}{\indent\sc Remark}

\newcommand{\skal}[2]{\langle #1,#2\rangle}

\begin{document}

\title{Two notes on harmonic distributions}

\author{Kamil Niedzia\l omski}

\subjclass[2000]{53C43, 58E20}

\keywords{Harmonic section, distribution, tangent bundle, Sasaki metric, conformal deformation}

\thanks{The author was supported by the Polish NSC grant No. 6065/B/H03/2011/40}

\address{
Department of Mathematics and Computer Science \endgraf
University of \L\'{o}d\'{z} \endgraf
ul. Banacha 22, 90-238 \L\'{o}d\'{z} \endgraf
Poland
}

\email{kamiln@math.uni.lodz.pl}

\begin{abstract}
We say that a distribution is harmonic if it is harmonic when considered as a section of a Grassmann bundle. We find new examples of harmonic distributions and show nonexistense of harmonic distrubutions on some Riemannian manifolds by two different approaches. Firstly, we lift distributions to the second tangent bundle equipped with the Sasaki metric. Secondly, we deform conformally the metric on a base manifold.  
\end{abstract}

\maketitle

\section{Introduction}

Harmonic map $\sigma:M\to N$ between Riemannian manifolds is a critical point of energy functional
\begin{equation*}
E(\sigma)=\int_M |\sigma_{\ast}|^2\,{\rm vol}_M.
\end{equation*}
When the considered map $\sigma$ is a section of a submersion, we may define the weaker condition of harmonicity. In the tangent space of any submersion we may distinguish the vertical subspace and therefore we may consider the vertical projection $\sigma_{\ast}^V(X)$ of a vector $\sigma_{\ast}(X)$, $X\in T_xM$. Thus, we define harmonicity of a section via vanishing of the Euler--Lagrange equation of the vertical energy functional
\begin{equation*}
E^V(\sigma)=\int_M |\sigma_{\ast}^V|^2\,{\rm vol}_M.
\end{equation*}

Harmonic distribution (plane field) $\sigma$ on a Riemannian manifold $M$ is a distribution, which considered as a section $\sigma:M\to {\rm Gr}_p(M)$, $p=\dim \sigma$, of a Grassmann bundle is harmonic (We equip ${\rm Gr}_p(M)$ with the Riemannian metric induced by the Riemannian metric on $M$ and the invariant inner product on $O(n)$). The following distributions are known to be harmonic:
\begin{enumerate}
\item\label{hd1} generalized Hopf fibrations \cite{cy,ggv},
\item\label{hd2} characteristic distribution of a contact structure \cite{ggv},
\item\label{hd3} invariant distributions on compact semisimple Lie groups \cite{ch}.
\item\label{hd4} distributions orthogonal to harmonic unit vector fields \cite{ggv}.
\end{enumerate}
Notice that example \eqref{hd2} is a special case of example \eqref{hd4} since Reeb vector field of a contact structure is unit harmonic. 

In this paper we find new examples of harmonic distributions and show nonexistence of harmonic distributions on some Riemannian manifolds. We consider two different approaches:

Firstly, we define two types of lifts of a distribution to the second tangent bundle: vertical and harmonic lift. We equip second tangent bundle with the Sasaki metric. Using the formula for the Levi--Civita connection of Sasaki metric, we obtain the formula for harmonicity of a lift. By this formula we get that harmonicity of a lift implies harmonicity of a base distribution. Moreover, we show that in the case of constant sectional curvature the vertical lift is never harmonic. Moreover, for the vertical and horizontal lift of the tangent bundle, harmonicity is equivalent to the vanishing of the divergence of curvature tensor.  

Secondly, we modify conformally the Riemannian metric on a base manifold. We derive a condition for harmonicity with respect to the metric in conformal class. We show that in for a foliation by curves on a surface harmonicity does not depend on a choice of a conformal metrics. We give relevant examples. We consider the well konwn example of harmonic distribution -- the Hopf fibration.

Throughout the paper we will use the following index convention
\[
1\leq\alpha,\beta,\gamma\leq n,\quad 1\leq a,b,c\leq p,\quad p+1\leq i,j,k\leq n.
\]

\section{Harmonicity a distribution}

Let us begin with the definition of harmonicity of arbitrary section. Let $\pi:P\to M$ be a submersion. Assume $P$ and $M$ are equipped with Riemannian metrics and let $\nabla^P$ and $\nabla^M$ denote the corresponding Levi--Civita connections, respectively. The {\it vertical distribution} $\mathcal{V}\subset TP$ is the kernel of a differential of the projection,
\begin{equation*}
\mathcal{V}=\ker \pi_{\ast}.
\end{equation*}
If $\sigma:M\to P$ is a section, then the $\mathcal{V}$--component of a vector $\sigma_{\ast}(X)$, $X\in T_xM$, is denoted by $\sigma_{\ast}^V(X)\in \mathcal{V}_{\sigma(x)}$. In the bundle $T^*M\times\sigma(\mathcal{V})$, where $\sigma(\mathcal{V})\to M$ is a pull-back bundle, define the metric induced from the Riemannian metrics on $M$ and $P$. We define {\it vertical energy functional} as follows
\begin{equation*}
E^V(\sigma)=\frac{1}{2}\int_M |\sigma_{\ast}^V|^2\,{\rm vol}_M
\end{equation*}
if $M$ is compact, otherwise we integrate over all compact sets $U\subset M$ to define $E^V(\sigma,U)$. We say that $\sigma$ is {\it harmonic} if the first variation of the vertical energy functional is equal to zero.

In the pull--back bundle $\sigma^{-1}(TP)$ there is only one connection $\nabla^{\sigma}$ such that
\begin{equation*}
\nabla^{\sigma}_X(Y\circ\sigma)=\nabla^P_{\sigma_{\ast}X}Y,\quad X\in T_xM,\quad Y\in\Gamma(TP).
\end{equation*} 
(see \cite{bw}). Then connections $\nabla^M$ and $\nabla^{\sigma}$ induce connection $\nabla$ in the bundle $T^*M\times\sigma^{-1}(TP)$. Its restriction to $T^*M\times \sigma^{-1}(\mathcal{V})$ is denoted by $\nabla^V$. Notice, that if the fibers $\mathcal{V}_x$, $x\in M$, are totally geodesic, then $\nabla^V$ coincides with $\nabla$. We define the {\it vertical tension field} $\tau^V(\sigma)\in\Gamma(\sigma^{-1}(\mathcal{V}))$ of a section $\sigma$ as follows
\begin{equation*}
\tau^V(\sigma)={\rm trace}\nabla^V\sigma_{\ast}.
\end{equation*} 
It appears \cite{wo} that the section is harmonic if and only if its vertical tension field vanishes.

Now we turn to harmonicity of a distribution. Let $\sigma$ be a $k$--dimensional distribution on a $n$--dimensional Riemannian manifold $M$. Let $G_k(M)_x$ be the space of all $k$--dimensional subspaces of the tangent space $T_xM$. Then $G_k(M)$ are the fibers if the {\it Grassmann bundle} $\pi:G_k(M)\to M$. Let $G=O(n)$ and $H=O(k)\times O(n-k)\subset G$. Then, each fiber $G_k(M)$ is isomorphic to the quotient $G/H$. Hence the Grassmann bundle is the associated bundle $G_k(M)=O(M)\times_G (G/H)$, where $O(M)$ is the principal bundle of all orthonormal frames on $M$. Let $\mathfrak{g}$ and $\mathfrak{h}$ denote the Lie algebras of $G$ and $H$ respectively. The subspace
\begin{equation*}
\mathfrak{m}=\left\{\left(\begin{array}{cc} 0 & A \\ A^{\top} & 0 \end{array}\right)\mid \textrm{$A$ is any $(n-k)\times k$ matrix}\right\}
\end{equation*}
is the complement of $\mathfrak{h}$ in $\mathfrak{g}$ i.e. $\mathfrak{g}=\mathfrak{h}\oplus\mathfrak{m}$. Let $\skal{\cdot}{\cdot}$ denotes $G$--invariant metric on $G/H$ or equivalently ${\rm Ad}_G(H)$--invariant inner product on $\mathfrak{m}$,
\begin{equation*}
\skal{A}{B}=-\frac{1}{2}{\rm trace}(AB),\quad A,B\in\mathfrak{m}.
\end{equation*}
The vertical subspace $\mathcal{V}_x\subset TG_k(M)$, $x\in M$, is isomorphic to $\mathfrak{m}$, hence the inner product $\skal{\cdot}{\cdot}$ induces an inner product in $\mathcal{V}$. 

Distribution $\sigma$ can be seen as a map $\sigma:M\to G_k(M)$, hence as a section of the Grassmann bundle. Thus we may speak about harmonicity of $\sigma$. 

Now we want to express the vertical tension field $\tau^V(\sigma)$ (for more details see \cite{wo} and \cite{cy}). Let $\pi_{O(\sigma)}:O(\sigma)\to M$ be a subbundle of $\pi_{O(M)}:O(M)\to M$ of orthonormal frames adopted to $\sigma$. Let
\begin{equation*}
\rho:O(M)\to O(M)/H=G_k(M)
\end{equation*}  
be a projection. Let ${\rm Ad}_G$ denotes the adjoint representation of $G$ on $\mathfrak{g}$. Consider the following bundles over $M$
\begin{align*}
&\mathfrak{G}=O(M)\times_{{\rm Ad}_G}\mathfrak{g},\\
&\mathfrak{H}=O(\sigma)\times_{{\rm Ad}_G}\mathfrak{h},\\
&\mathfrak{M}=O(\sigma)\times_{{\rm Ad}_G}\mathfrak{m}.
\end{align*}
Each of above bundles can be identified with the subbundle of ${\rm End}(TM)$ via the map
\begin{equation*}
u\cdot A\mapsto u A u^{-1}:T_xM\to T_xM,
\end{equation*}
where the frame $u\in O(M)$ or $u\in O(\sigma)$ is understood as a map $u:\mathbb{R}^n\to T_xM$. With respect to this identification, $\mathfrak{H}$ is the bundle of all skew--symmetric endomorphism preserving $\sigma$ and $\sigma^{\bot}$, whereas $\mathfrak{M}$ the bundle of all skew--symmetric endomorphisms of the form $(L,L^{\top})$, where $L:\sigma\to \sigma^{\bot}$ and $L^{\top}:\sigma^{\bot}\to\sigma$ is the adjoint to $L$. Moreover
\begin{equation*}
\pi_{O(M)}^{-1}\mathfrak{G}\equiv \mathcal{V}^{O(M)},\quad \pi_{O(\sigma)}^{-1}\mathfrak{H}\equiv \mathcal{V}^{O(\sigma)} 
\end{equation*} 
via the identification
\begin{equation*}
u\cdot A \mapsto A^*_u,
\end{equation*}
where $A^*_u$ is the fundamental vertical vector. Moreover
\begin{equation*}
\pi_{O(\sigma)}^{-1}\mathfrak{H}\oplus\pi_{O(\sigma)}^{-1}\mathfrak{M}=\pi_{O(\sigma)}^{-1}\mathfrak{G}
=\mathcal{V}^{O(M)}|O(\sigma).
\end{equation*}
Since $\rho|O(\sigma):O(\sigma)\to G_k(M)$ is constant, then $\ker\rho_{\ast}=\mathcal{V}^{O(\sigma)}$. Since, by above,
\begin{equation*}
\mathcal{V}^{O(M)}|O(\sigma)=\mathcal{V}^{O(\sigma)}\oplus\pi_{O(\sigma)}^{-1}\mathfrak{M},
\end{equation*} 
we get
\begin{equation*}
\mathcal{V}\equiv \pi^{-1}\mathfrak{M}\quad\textrm{and}\quad \sigma^{-1}(\mathcal{V})\equiv \mathfrak{M}.
\end{equation*}
Finally, we may identify the pull-back of the vertical distribution $\sigma^{-1}(\mathcal{V})$ with the bundle of skew--symmetric endomorphisms $(L,L^{\top})$ of $TM$ described above. One can show \cite{wo} that element $\sigma_{\ast}^V(X)$, $X\in T_xM$, is identified with the endomorphism
\begin{equation*}
P_X(Y)=(\nabla^M_XY^{\bot})^{\top}+(\nabla^M_XY^{\top})^{\bot},
\end{equation*}
where $\top$ and $\bot$ denote projection onto $\sigma$ and $\sigma^{\bot}$, respectively. Hence, the vector space $T^*_xM\otimes \sigma^{-1}(\mathcal{V})_x$ is identified with the space of $\mathbb{R}$--bilinear maps $T_xM\times T_xM\ni(X,Y)\mapsto P_X(Y)$ and the connection is identified with the connection induced from the Levi--Civita connection $\nabla^M$. Thus, the vertical tension field $\tau^V(\sigma)$ is equal \cite{cy}
\begin{equation*}
\tau^V(\sigma)={\rm trace}\nabla^V \sigma_{\ast}=\sum_{\alpha}\nabla^M_{e_{\alpha}}P_{e_{\alpha}}-P_{\nabla^M_{e_{\alpha}}e_{\alpha}},
\end{equation*}
where $e_1,\ldots,e_n$ is a local orthonormal basis on $M$. Since $P$ is determined by its values on vectors tangent to $\sigma$, it suffices to evaluate $\tau^V(\sigma)$ such vectors. For $X\in\sigma$ we have
\begin{equation}\label{eq:vtf}
\begin{split}
\tau^V(\sigma)(X) &=\sum_{\alpha} (\nabla_{e_{\alpha}}P_{e_{\alpha}})(X)-P_{\nabla^M_{e_{\alpha}}e_{\alpha}}(X)\\
&=\sum_{\alpha}\nabla^M_{e_{\alpha}}(P_{e_{\alpha}}(X))-P_{e_{\alpha}}(\nabla^M_{e_{\alpha}}X)-P_{\nabla_{e_{\alpha}}e_{\alpha}}(X)\\
&=\sum_{\alpha} \nabla^M_{e_{\alpha}}(\nabla^M_{e_{\alpha}}X)^{\bot}-(\nabla_{e_{\alpha}}(\nabla_{e_{\alpha}}X)^{\top})^{\bot}
-(\nabla_{\nabla_{e_{\alpha}}e_{\alpha}}X)^{\bot}\\
&-\sum_{\alpha}(\nabla_{e_{\alpha}}(\nabla_{e_{\alpha}}X)^{\bot})^{\top}\\
&=\sum_{\alpha} (\nabla^M_{e_{\alpha}}(\nabla^M_{e_{\alpha}}X)^{\bot})^{\bot}-(\nabla_{e_{\alpha}}(\nabla_{e_{\alpha}}X)^{\top})^{\bot}
-(\nabla_{\nabla_{e_{\alpha}}e_{\alpha}}X)^{\bot}\\
&=\sum_{\alpha}(\nabla_{e_{\alpha}}\nabla_{e_{\alpha}}X-2\nabla_{e_{\alpha}}(\nabla_{e_{\alpha}}X)^{\top}-\nabla_{\nabla_{e_{\alpha}}e_{\alpha}}X)^{\bot} \\
&=\sum_{\alpha}((\nabla^2X)(e_{\alpha},e_{\alpha})-2\nabla_{e_{\alpha}}(\nabla_{e_{\alpha}}X)^{\top})^{\bot}.
\end{split}
\end{equation}
We put
\begin{align*}
\tau^V(\sigma)(X,Y)=g(\tau^V(\sigma)(X),Y),\quad X\in\sigma,\quad Y\in\sigma^{\bot}.
\end{align*}
where $g$ is the Riemannian metric on $M$. Notice that $\tau^V(\sigma)(X,Y)$ is tensorial with respect to $X$ and $Y$.

\begin{proposition}\label{pr:vertenfieldxy}
The vertical tension field $\tau^V(\sigma)(X,Y)$ for $X\in \sigma$ and $Y\in \sigma^{\bot}$ is equal to
\begin{equation}\label{eq:vertenfieldxy}
\begin{split}
\tau^V(\sigma)(X,Y) &=\sum_bg((\nabla^2X)(e_b,e_b),Y)-\sum_jg((\nabla^2Y)(e_j,e_j),X) \\
&+2\sum_{b,c}g(X,\nabla_{e_b}e_c)g(\nabla_{e_b}e_c,Y)-2\sum_{j,k}g(X,\nabla_{e_j}e_k)g(\nabla_{e_j}e_k,Y).
\end{split}
\end{equation}
\end{proposition}
\begin{proof}
For any $X\in\sigma$ and $Y\in \sigma^{\bot}$ by \eqref{eq:vtf}
\begin{align*}
\tau^V(\sigma)(X,Y)&=\sum_{\alpha}g((\nabla^2X)(e_{\alpha},e_{\alpha})-2\nabla_{e_{\alpha}}(\nabla_{e_{\alpha}}X)^{\top},Y) \\
&=\sum_bg((\nabla^2X)(e_b,e_b),Y)+\sum_jg((\nabla^2X)(e_j,e_j),Y) \\
&+2\sum_bg((\nabla_{e_b}X)^{\top},\nabla_{e_b}Y)+2\sum_jg((\nabla_{e_j}X)^{\top},\nabla_{e_j}Y).
\end{align*}
Since
\[
g(\nabla_{e_j}\nabla_{e_j}X,Y)=-2g(\nabla_{e_j}X,\nabla_{e_j}Y)-g(X,\nabla_{e_j}\nabla_{e_j}Y),
\]
then 
\[
g((\nabla^2X)(e_j,e_j),Y)=-g((\nabla^2Y)(e_j,e_j),X)-2g(\nabla_{e_j}X,\nabla_{e_j}Y).
\]
Therefore
\begin{align*}
\tau^V(\sigma)(X,Y) &=\sum_bg((\nabla^2X)(e_b,e_b),Y)-\sum_jg((\nabla^2Y)(e_j,e_j),X) \\
&-2\sum_jg(\nabla_{e_j}X,\nabla_{e_j}Y)+2\sum_bg((\nabla_{e_b}X)^{\top},\nabla_{e_b}Y) \\
&+2\sum_jg(\nabla_{e_j}e_a,\nabla_{e_j}e_i)-2\sum_jg((\nabla_{e_j}e_a)^{\bot},\nabla_{e_j}e_i) \\
&=\sum_bg((\nabla^2X)(e_b,e_b),Y)-\sum_jg((\nabla^2Y)(e_j,e_j),X) \\
&+2\sum_{b,c}g(\nabla_{e_b}X,e_c)g(e_c,\nabla_{e_b}Y)-2\sum_{j,k}g(\nabla_{e_j}X,e_k)g(e_k,\nabla_{e_j}Y) \\
&=\sum_bg((\nabla^2X)(e_b,e_b),Y)-\sum_jg((\nabla^2Y)(e_j,e_j),X) \\
&+2\sum_{b,c}g(X,\nabla_{e_b}e_c)g(\nabla_{e_b}e_c,Y)-2\sum_{i,k}g(X,\nabla_{e_j}e_k)g(\nabla_{e_j}e_k,Y),
\end{align*}
which ends the proof.
\end{proof}

Since the vanishing of the vertical tension field is equivalent to harmonicity of a section, we have that a distribution $\sigma$ is harmonic if and only if
\begin{equation*}
\tau^V(\sigma)(X,Y)=0\quad\textrm{for all $X\in \sigma$, $Y\in \sigma^{\bot}$}.
\end{equation*} 

By the formula \ref{eq:vertenfieldxy} for the vertical tension field we immediately get the following result \cite{ggv}.

\begin{proposition}\label{symverticaltension}
For $X\in\sigma$ and $Y\in\sigma^{\bot}$ we have
\begin{equation*}
\tau^V(\sigma)(X,Y)=-\tau^V(\sigma^{\bot})(Y,X).
\end{equation*}
In particular, $\sigma$ is harmonic if and only if $\sigma^{\bot}$ is harmonic.
\end{proposition} 

\section{Lifts of distribution to the tangent bundle}
In this section, we construct the lifts of a distribution to tangent bundle and compute the tension field of obtained distributions. By this procedure, we obtain new examples of harmonic distributions.

Let $TTM$ be the second tangent bundle. The Levi--Civita connection $\nabla^M$ and projection $\pi_{TM}:TM\to M$ induce the following decomposition \cite{ko}
\begin{equation*}
T_{\xi}TM=\mathcal{H}^{TM}_{\xi}\oplus\mathcal{V}^{TM}_{\xi},\quad \xi\in T_xM.
\end{equation*} 
Then, for any $X\in T_xM$ there is unique horizontal lift $X^h_{\xi}\in\mathcal{H}^{TM}_{\xi}$ and unique vertical lift $X^v_{\xi}\in\mathcal{V}^{TM}_{\xi}$. Let $g_S$ denotes the Sasaki metric on $TM$,
\begin{align*}
g_S(X^h,Y^h) &=g(X,Y), \\
g_S(X^h,Y^v) &=0, \\
g_S(X^v,Y^v) &=g(X,Y).
\end{align*}
The corresponding Levi--Civita connection $\nabla^{TM}$ is of the form \cite{gk}
\begin{align*}
\left(\nabla^{TM}_{X^h}Y^h\right)_{\xi} &=(\nabla^M_XY)^h_{\xi}-\frac{1}{2}(R(X,Y)\xi)^v_{\xi},\\
\left(\nabla^{TM}_{X^h}Y^v\right)_{\xi} &=(\nabla^M_XY)^v_{\xi}+\frac{1}{2}(R(\xi,Y)X)^h_{\xi},\\
\left(\nabla^{TM}_{X^v}Y^h\right)_{\xi} &=\frac{1}{2}(R(\xi,X)Y)^h_{\xi},\\
\left(\nabla^{TM}_{X^v}Y^v\right)_{\xi} &=0.
\end{align*} 

Let $\sigma$ be a $p$--dimensional distribution on $M$. Put
\begin{align*}
&\sigma^h(\xi)={\rm span}\{(e_1)^h_{\xi},\ldots,(e_p)^h_{\xi}\}, \\
&\sigma^v(\xi)={\rm span}\{(e_1)^v_{\xi},\ldots,(e_p)^v_{\xi}\},
\end{align*}
where $\sigma(x)={\rm span}\{e_1,\ldots,e_p\}$, $\pi_{TM}(\xi)=x$. Then $\sigma^h$ and $\sigma^v$ are two $p$--dimensional distributions on $TM$ called {\it horizontal} and {\it vertical lift} of $\sigma$, respectively.

Recall, that the divergence of a tensor field $T$ of type $(s+1,1)$ is a tensor field ${\rm div}T$ of type $(s,1)$ of the form
\begin{equation*}
{\rm div}T(X_1,\ldots,X_s)=\sum_{\alpha}(\nabla_{e_{\alpha}}T)(e_{\alpha},X_1,\ldots,X_s),
\end{equation*}
where $e_1,\ldots,e_n$ is a local orthonormal basis.

\begin{theorem}\label{liftdistributiontension}
The vertical tension field of lifts of distribution $\sigma$ are equal
\begin{enumerate}
\item for the horizontal lift $\sigma^h$:
\begin{equation}\label{eq:h1}
\tau^V(\sigma^h)(X^h_{\xi},Y^v_{\xi})=-\frac{1}{2}g(({\rm div}R)(X,\xi),Y)-\sum_{\alpha}g(R(e_{\alpha},(\nabla^M_{e_{\alpha}}X)^{\bot})\xi,Y), \end{equation}
where $X\in\sigma$, $Y\in TM$, and
\begin{equation} \label{eq:h2}
\tau^V(\sigma^h)(X^h_{\xi},Y^h_{\xi})=\tau^V(\sigma)(X,Y)-\sum_i g(R(e_i,X)\xi,R(e_i,Y)\xi),
\end{equation}
where $X\in\sigma$, $Y\in\sigma^{\bot}$,
\item for the vertical lift $\sigma^v$:
\begin{equation}\label{eq:v1}
\tau^V(\sigma^v)(X^v_{\xi},Y^h_{\xi})=\frac{1}{2}g(({\rm div}R)(Y,\xi),X)+\sum_{\alpha}g(R(\xi,(\nabla^M_{e_{\alpha}}X)^{\bot})e_{\alpha},Y),
\end{equation}
where $X\in\sigma$, $Y\in TM$, and
\begin{equation}\label{eq:v2}
\tau^V(\sigma^v)(X^v_{\xi},Y^v_{\xi})=\tau^V(\sigma)(X,Y)-\frac{1}{4}\sum_{\alpha}g(R(\xi,X)e_{\alpha},R(\xi,Y)e_{\alpha}),
\end{equation}
where $X\in\sigma$, $Y\in\sigma^{\bot}$.
\end{enumerate}
In particular,
\begin{equation*}
\tau^V(\mathcal{H}^{TM})(X^h_{\xi},Y^v_{\xi})=-\tau^V(\mathcal{V}^{TM})(Y^v_{\xi},X^h_{\xi})=-\frac{1}{2}g(({\rm div}R)(X,\xi),Y).
\end{equation*}
\end{theorem} 
\begin{proof}
{\it Proof of \eqref{eq:h1}}. Extend $X$ any $Y$ to local vector fields such that $g(X,e_a)=0$ and $g(Y,e_{\alpha})=0$. Now, it suffices to use the formulae for $\nabla^{TM}$. 

{\it Proof of \eqref{eq:h2}}. By a direct computation, using the formulae for $\nabla^{TM}$, we get
\begin{align*}
\tau^V(\sigma^h)(X^h_{\xi},Y^v_{\xi}) &=\tau^V(\sigma)(X,Y)-\frac{1}{4}\sum_a g(R(e_a,X)\xi,R(e_a,Y)\xi)\\
&+\frac{1}{4}\sum_{\alpha}g(R(\xi,e_{\alpha})X,R(\xi,e_{\alpha})Y)
-\frac{3}{4}\sum_i g(R(e_i,X)\xi,R(e_i,Y)\xi).
\end{align*}
Now, \eqref{eq:h2} follows from the equality
\begin{align*}
\sum_{\alpha}g(R(\xi,e_{\alpha})X,R(\xi,e_{\alpha})Y) &=\sum_{\alpha,\beta}g(R(\xi,e_{\alpha})X,e_{\beta})g(R(\xi,e_{\alpha})Y,e_{\beta})\\
&=\sum_{\alpha,\beta}g(R(X,e_{\beta})\xi,e_{\alpha})g(R(Y,e_{\beta})\xi,e_{\alpha})\\
&=\sum_{\beta}g(R(X,e_{\beta})\xi,R(Y,e_{\beta})\xi).
\end{align*}

{\it Proof of \eqref{eq:v1}}. It suffices to extend $Y$ to local vector field such that $g(Y,e_{\alpha})=0$ and use the formula for the Levi--Civita connection $\nabla^{TM}$.

{\it Proof of \eqref{eq:v2}}. We have
\begin{align*}
\tau^V(\sigma^v)(X^v,Y^v) &=-\sum_{\alpha} g((\nabla^M)^2 Y,X)-2\sum_{\alpha,i}g(X,\nabla^M_{e_{\alpha}}e_i)g(Y,\nabla^M_{e_{\alpha}}e_i)\\
&-\frac{1}{4}\sum_{\alpha}g(R(\xi,X)e_{\alpha},R(\xi,Y)e_{\alpha}).
\end{align*}
Extending $Y$ to a local vector field such that $g(Y,e_i)=0$, we get
\begin{equation*}
\sum_{\alpha,i}g(X,\nabla^M_{e_{\alpha}}e_i)g(Y,\nabla^M_{e_{\alpha}}e_i)=
\sum_{\alpha}g((\nabla^M_{e_{\alpha}}X)^{\bot},(\nabla^M_{e_{\alpha}}Y)^{\bot}).
\end{equation*}
By the formula \eqref{eq:vtf} we get
\begin{equation*}
\tau^V(\sigma^{\bot})(Y,X)=\sum_{\alpha}g((\nabla^M)^2Y,X)+
2\sum_{\alpha}g((\nabla^M_{e_{\alpha}}X)^{\bot},(\nabla^M_{e_{\alpha}}Y)^{\bot}).
\end{equation*}
Combining above equalities and by proposition \ref{symverticaltension} we get \eqref{eq:v2}.
\end{proof}

As an immediate consequence of above theorem we get.
\begin{corollary}
Distributions $\mathcal{H}^{TM}$ and $\mathcal{V}^{TM}$ are harmonic if and only if the curvature tensor $R$ on $M$ is divergence free. In particular, $\mathcal{H}^{TM}$ and $\mathcal{V}^{TM}$ are harmonic on the symmetric space and on a manifold of constant scalar curvature. 
\end{corollary}

Moreover, the harmonicity of a lift implies harmonicity of a base distribution. 

\begin{corollary}\label{liftconstant}
Let $\sigma\neq TM$. Then
\begin{enumerate}
\item If $\sigma^h$ is harmonic, then $\sigma$ is harmonic.
\item If $\sigma^v$ is harmonic, then $\sigma$ is harmonic and $R(\sigma,\sigma^{\bot})=0$. In particular, on a manifold of nonzero constant sectional curvature the vertical lift $\sigma^v$ is never harmonic.
\end{enumerate}
\end{corollary}
\begin{proof}
Assume first $\sigma^h$ is harmonic. Replacing $\xi$ by $t\xi$, condition \eqref{eq:h2} is equivalent to the vanishing of a polynomial
\begin{equation*}
\chi(t)=\tau^V(\sigma)(X,Y)-t^2\sum_i g(R(e_i,X)\xi,R(e_i,Y)\xi).
\end{equation*}  
Hence $\tau^V(\sigma)(X,Y)=0$ for $X\in\sigma$ and $Y\in\sigma^{\bot}$. Therefore $\sigma$ is harmonic.

Assume now $\sigma^v$ is harmonic. Replacing $\xi$ by $t\xi$ in the condition \eqref{eq:v2} the vanishing of obtained polynomial implies
\begin{equation*}
\tau(\sigma)(X,Y)=0\quad\textrm{and}\quad\sum_{\alpha}g(R(\xi,X)e_{\alpha},R(\xi,Y)e_{\alpha})=0
\end{equation*}
for $X\in\sigma$ and $Y\in\sigma^{\bot}$. Hence $\sigma$ is harmonic. Moreover, putting $\xi=X+Y$ in the second condition, we get $|R(X,Y)e_{\alpha}|^2=0$. Thus $R(\sigma,\sigma^{\bot})=0$. Since the curvature tensor $R$ on a manifold of nonzero constant sectional curvature $\kappa$ is equal $R(X,Y)Z=\kappa(g(Z,Y)X-g(Z,Y)X)$, it follows that in this case $R(\sigma,\sigma^{\bot}\neq 0$. 
\end{proof}

\section{Conformal deformations}

In this section we derive the formula for the tension field under conformal deformation of a Riemannian metric.  

Let $(M,g)$ be a Riemannian manifold. Consider a Riemannian metric $\tilde g=e^{2\mu}g$, where $\mu$ a smooth function on $M$. Let $\nabla$ and $\tilde\nabla$ denote the Levi--Civita connections of $g$ and $\tilde g$, respectively. For any $X,Y,Z\in\Gamma(TM)$ one may check that $\nabla$ and $\tilde\nabla$ are related as follows
\begin{equation}\label{confnabla}
\tilde\nabla_XY=\nabla_XY+(Y\mu)X+(X\mu)Y-g(X,Y)\nabla\mu
\end{equation}

Let $\sigma$ be a $p$--dimensional distribution on $M$ and $\sigma^{\bot}$ denotes its orthogonal complement. Put $q=\dim\sigma^{\bot}$. The {\it mean curvature} of $\sigma$ with respect to $g$ is a vector field $H_{\sigma}$ of the form
\begin{equation*}
H_{\sigma}=\sum_a \left(\nabla_{e_a}e_a\right)^{\bot},
\end{equation*}
where $e_1,\ldots,e_p$ is an orthonormal basis of $\sigma$ (with respect to $g$).

Denote the vertical tension fields of $\sigma$ with respect to $g$ and $\tilde g$ by $\tau^V(\sigma)$ and $\tilde\tau^V(\sigma)$, respectively. 

\begin{theorem}
Vertical tension fields $\tau^V(\sigma)$ and $\tilde\tau^V(\sigma)$ are related as follows
\begin{gather}\label{harmconf2}
\begin{split}
e^{2\mu}\tilde\tau^V(\sigma)(X,Y) &= \tau^v(\sigma)(X,Y)+(p-q)g(\nabla\mu,X)g(\nabla\mu,Y) \\
&-2g(\nabla\mu,X)g(H_{\sigma},Y)+2g(\nabla\mu,Y)g(H_{\sigma^{\bot}},X) \\
&-2g(\nabla_Y(\nabla\mu)^{\bot},X)+2g(\nabla_X(\nabla\mu)^{\top},Y) \\
&+(n-2)g(\nabla_{\nabla\mu}X,Y),
\end{split}
\end{gather} 
where $X\in \sigma$ and $Y\in\sigma^{\bot}$. 
\end{theorem}
\begin{proof}
Let $e_1,\ldots,e_n$ be a local orthonormal basis with respect to $g$ such that $\sigma_a\in\sigma$. Put $\mu_{\alpha}=g(\nabla\mu,e_{\alpha})$. The basis $f_{\alpha}=e^{-\mu}e_{\alpha}$ is orthonormal with respect to $\tilde g$. Since the condition $\tilde\tau^V(\sigma)(X,Y)$ is tensorial with respect to $X\in\sigma$ and $Y\sigma^{\bot}$ is suffices to prove the formula \eqref{harmconf2} for $X=f_a$ and $Y=f_i$. We have
\begin{align}
\tilde g(\tilde{\nabla}_{f_b}f_c,f_a) &=e^{-\mu}(g(\nabla_{e_b}e_c,e_a)+\mu_c\delta_{ab}-\mu_a\delta_{bc}), \label{mth1} \\
\tilde g(\tilde{\nabla}_{f_b}f_c,f_i) &=e^{-\mu}(g(\nabla_{e_b}e_c,e_i)-\mu_i\delta_{bc}), \label{mth2}
\end{align}
Put
\begin{align*}
P_1 &=2e^{2\mu}\sum_{b,c}\tilde g(\tilde\nabla_{f_b}f_c,f_a)\tilde g(\tilde\nabla_{f_b}f_c,f_i), \\
P_2 &=2e^{2\mu}\sum_{j,k}\tilde g(\tilde\nabla_{f_j}f_k,f_i)\tilde g(\tilde\nabla_{f_j}f_k,f_a),
\end{align*}
and
\begin{align*}
Q_1 &=e^{2\mu}\sum_b\tilde g(\tilde\nabla_{f_b}\tilde\nabla_{f_b}f_a,f_i), \\
Q_2 &=e^{2\mu}\sum_j\tilde g(\tilde\nabla_{f_j}\tilde\nabla_{f_j}f_i,f_a),
\end{align*}
and
\begin{align*}
S_1 &=e^{2\mu}\sum_b\tilde g(\tilde\nabla_{\tilde\nabla_{f_b}f_b}f_a,f_i), \\
S_2 &=e^{2\mu}\sum_j\tilde g(\tilde\nabla_{\tilde\nabla_{f_j}f_j}f_i,f_a),
\end{align*}
Then, by \eqref{confnabla} using \eqref{mth1} and \eqref{mth2} we get
\begin{gather*}
\begin{split}
P_1 &=\sum_{b,c}g(e_a,\nabla_{e_b}e_c)g(\nabla_{e_b}e_c,e_i)-\mu_i g(\sum_b\nabla_{e_b}e_b,e_a) \\
&+g(\nabla_{e_a}(\nabla\mu)^{\top},e_i)-\mu_ag(H_{\sigma},e_i)+(p-1)\mu_a\mu_i
\end{split}
\end{gather*}
and
\begin{gather*}
\begin{split}
Q_1 &=\sum_b g(\nabla_{e_b}\nabla_{e_b}e_a,e_i)+(2-p)\mu_a\mu_i-g(\nabla_{(\nabla\mu)^{\top}}e_a,e_i) \\
&+\mu_i g(\sum_b\nabla_{e_b}e_b,e_a)+\mu_ag(H_{\sigma},e_i)-{\rm hess}_{\mu}(e_a,e_i)
\end{split}
\end{gather*}
and 
\begin{gather*}
\begin{split}
S_1 &=\sum_b g(\nabla_{\nabla_{e_b}e_b}e_a,e_i)+\mu_a g(H_{\sigma},e_i)-\mu_a\mu_i \\
&-\mu_i g(\sum_b\nabla_{e_b}e_b,e_a)+g(\nabla_{(\nabla\mu)^{\top}}e_a,e_i)-pg(\nabla_{\nabla\mu}e_a,e_i),
\end{split}
\end{gather*}
where ${\rm hess}_{\mu}(e_a,e_i)=g(\nabla_{e_a}\nabla\mu,e_i)$ denotes the hessian of $\mu$ in the direction of $e_a$ and $e_i$.

Analogously, interchanging $i$ with $a$, $b$ with $j$ and $\top$ with $\bot$, we get
\begin{gather*}
\begin{split}
P_2 &=\sum_{j,k}g(e_a,\nabla_{e_j}e_k)g(\nabla_{e_j}e_k,e_i)-\mu_ag(\sum_j\nabla_{e_j}e_j,e_i) \\
&+g(\nabla_{e_i}(\nabla\mu)^{\bot},e_a)-\mu_ig(H_{\sigma^{\bot}},e_a)+(q-1)\mu_a\mu_i
\end{split}
\end{gather*}
and
\begin{gather*}
\begin{split}
Q_2 &=\sum_j g(\nabla_{e_j}\nabla_{e_j}e_i,e_a)+(2-q)\mu_a\mu_i-g(\nabla_{(\nabla\mu)^{\bot}}e_i,e_a) \\
&+\mu_a g(\sum_j\nabla_{e_j}e_j,e_i)+\mu_ig(H_{\sigma^{\bot}},e_a)-{\rm hess}_{\mu}(e_i,e_a)
\end{split}
\end{gather*}
and 
\begin{gather*}
\begin{split}
S_2 &=\sum_j g(\nabla_{\nabla_{e_j}e_j}e_i,e_a)+\mu_i g(H_{\sigma^{\bot}},e_a)-\mu_a\mu_i \\
&-\mu_a g(\sum_j\nabla_{e_j}e_j,e_i)+g(\nabla_{(\nabla\mu)^{\bot}}e_i,e_a)-qg(\nabla_{\nabla\mu}e_i,e_a).
\end{split}
\end{gather*}
Since $e^{2\mu}\tau^v_{\tilde g}(\sigma)_{a,i}=P_1-P_2+(Q_1-S_1)+(Q_2-S_2)$, \eqref{harmconf2} holds.
\end{proof}

\begin{corollary}\label{harmconf2dim}
If $\sigma$ a foliation by curves on a $2$--dimensional manifold, then the harmonicity of $\sigma$ depends only on the conformal structure on $M$.
\end{corollary}
\begin{proof}
Let $\sigma={\rm Span}(X)$, $\sigma^{\bot}={\rm Span}(Y)$, where $X,Y$ is an orthonormal frame on $M$. Let $\mu_X=g(\nabla\mu,X)$, $\mu_Y=g(\nabla\mu,Y)$ and let $H_{\sigma}=h_{\sigma}Y$, $H_{\sigma^{\bot}}=h_{\sigma^{\bot}}X$. Then condition \eqref{harmconf2} simplifies to the following
\begin{gather*}
\begin{split}
e^{2\mu}\tilde\tau^v(\sigma)(X,Y) &=\tau^v(\sigma)(X,Y)-2\mu_Xh_{\sigma}+2\mu_Yh_{\sigma^{\bot}} \\
&-2g(\nabla_Y(\mu_YY),X)+2g(\nabla_X(\mu_XX),Y) \\
&=\tau^v(\sigma)(X,Y)-2\mu_Xh_{\sigma}+2\mu_Yh_{\sigma^{\bot}}\\
&-2\mu_Y g(\nabla_YY,X)+2\mu_X g(\nabla_XX,Y) \\
&=\tau^v(\sigma)(X,Y).
\end{split}
\end{gather*}
Thus condition $\tau^v(\sigma)=0$ is equivalent to $\tilde\tau^v(\sigma)=0$. Hence, harmonicity depends only on the conformal structure on $M$.
\end{proof}

\begin{remark0}
Corollary \ref{harmconf2dim} is analogous to the general fact that harmonicity of a map from $2$--dimensional manifold depends only on the conformal structure (see \cite[Corollary 3.5.4]{bw}).
\end{remark0}

\begin{corollary}\label{Cor:totgeod}
Let $\sigma$ be totally geodesic foliation on a Riemannian manifold $(M,g)$ and let $\mu$ be a function such that $(\nabla\mu)^{\bot}=0$. If $\sigma$ is harmonic then $\sigma$ is harmonic with respect to a metric $\tilde g=e^{\mu}g$.
\end{corollary}
\begin{proof}
Follows immediately by \eqref{harmconf2}.
\end{proof}

\begin{remark0}
Notice, that with the assumptions of above Corollary, the foliation $\sigma$ is also totally geodesic with respect to $\tilde g$.
\end{remark0}

In the end we show that formula \eqref{harmconf2} implies results on nonexistence of harmonic distributions. 

\begin{example}
Let $\sigma$ be a foliation by grate circles of the Hopf Fibration $S^1\to S^3\to S^2$. We may describe Hopf fibration as follows: Sphere $S^3\subset\mathbb{R}^4$ is a Lie group, whose Lie algebra is spanned by vectors $X,Y,Z$ such that 
\begin{equation}\label{eq:hopf}
[X,Y]=2Z,\quad [Y,Z]=2X,\quad [Z,X]=2Y.
\end{equation}
Then $\sigma$ is spanned by vector $X$. One can show that $\sigma$ is harmonic  with respect to standard inner product $g$ on $S^3$ \cite{ggv,cy}. Suppose there exists function $\mu$ such that $\sigma$ is harmonic with respect to $e^{2\mu}g$, where $\mu$ is not a constant. By Koszul formula for the Levi--Civita connection we have
\begin{equation*}
(\nabla_XX)^{\bot}=(\nabla_YY)^{\top}=(\nabla_ZZ)^{\top}=0,\quad g(\nabla_YX,Z)=g(\nabla_ZX,Y)=1.
\end{equation*}
Thus, by \eqref{eq:hopf}, we have
\begin{align*}
e^{2\mu}\tilde\tau^V(\sigma)(X,Y) &=-(X\mu)(Y\mu)-2(Z\mu)g(\nabla_YZ,X)+(Z\mu)g(\nabla_ZX,Y)\\
&=-(X\mu)(Y\mu)-Z\mu
\end{align*}
and
\begin{align*}
e^{2\mu}\tilde\tau^V(\sigma)(X,Z) &=-(X\mu)(Z\mu)-2(Y\mu)g(\nabla_ZY,X)+(Y\mu)g(\nabla_YX,Z)\\
&=-(X\mu)(Z\mu)+Y\mu.
\end{align*}
Hence $\tilde\tau^V=0$ if and only if $Y\mu=Z\mu=0$. Since $[Y,Z]=2X$, then $Y\mu=Z\mu=0$ implies $X\mu=0$. Thus $\mu$ is constant, which contradicts the assumption. Finally, there is no Riemannian metric, except for constant multiplicity of $g$, in the conformal class of $g$, such that $\sigma$ is harmonic.
\end{example}

\end{document}